\documentstyle{amsart}

\title[Existence results for nonlinear...]
{Existence results for nonlinear elliptic problems on fractal domains}

\author{Massimiliano Ferrara, Giovanni Molica Bisci and Du\v{s}an Repov\v{s}}

\address[M. Ferrara]{University of Reggio Calabria and CRIOS University Bocconi of Milan, Via dei Bianchi presso Palazzo Zani, 89127 Reggio Calabria, Italy}
\curraddr{}
\email{massimiliano.ferrara@@unirc.it}

\address[G. Molica Bisci]{Department P.A.U., Architecture Faculty, University of Reggio Calabria, 89124 - Reggio Calabria, Italy} \email{gmolica@@unirc.it}
\address[D. Repov\v{s}]{Faculty of Education, and Faculty of Mathematics and Physics - University of Ljubljana, POB 2964, Ljubljana, Slovenia 1001} \email{dusan.repovs@@guest.arnes.si}


\newtheorem{theorem}{Theorem}[section]
\newtheorem{lemma}{Lemma}[section]

\newtheorem{remark}{Remark}[section]
\newtheorem{example}{Example}[section]

\newcommand{\erre}{\mbox{\normalshape I\!R}}

\def\di{\displaystyle}

\def\R{{\rm I\!R}}
\def\RR{{\rm I\!R}}
\def\N{{\rm I\!N}}

\def\phi{\varphi}

\def\RR{{\rm I\!R}}
\def\N{{\rm I\!N}}

\def\phi{\varphi}

\def\phi{\varphi}

\def\di{\displaystyle}

\begin{document}


\begin{abstract}
Some existence results for a parametric Dirichlet problem defined on the Sierpi\'nski fractal are proved. More precisely, a critical point result for differentiable functionals is
exploited in order to prove the existence of a well determined open
interval of positive eigenvalues for which the problem admits at
least one non-trivial weak solution.
\end{abstract}

\maketitle


\section{Introduction}

\indent The purpose of the present paper is to establish some existence results for the following Dirichlet problem
\begin{equation}\tag{$S_{a,\lambda}^{f,g}$}\label{N}
\left\{
\begin{array}{l}
\Delta u(x)+a(x)u(x)=\lambda g(x)f(u(x)),\quad x \in V\setminus V_0, \\
u|_{V_0}=0,\\
\end{array}
\right.
\end{equation}
where $V$ stands for the Sierpi\'nski gasket in $(\erre^{N-1},|\cdot|)$, $N\geq 2$, $V_0$ is its intrinsic boundary (consisting of its $N$ corners), $\Delta$ denotes the weak Laplacian on $V$ and $\lambda$ is a positive real parameter. We assume that $f:\R\to\R$ is a continuous function and that the variable potentials $a, g:V\to \R$ satisfy the following conditions:
\begin{itemize}
\item[{$({\rm{h}}_1)$}] $a\in L^1(V,\mu)$ and $a\leq0$ almost everywhere in $V$;
\item[{$({\rm{h}}_2)$}] $g\in C(V)$ with $g\leq0$ and such that the restriction of $g$ to every open subset of $V$ is not identically zero.
\end{itemize}

Many physical problems on fractal domains lead to nonlinear models
(for example, reaction-diffusion equations, problems on elastic
fractal media or fluid flow through fractal regions), so it is
appropriate to study nonlinear partial differential equations on
fractals.\par
 In recent years there has been an increasing interest in
studying such equations, also motivated and stimulated by the
considerable amount of literature devoted to the definition of a
Laplacian operator for functions on fractal domains. Among the
contributions to the theory of nonlinear elliptic equations on
fractals we mention \cite{Fa99,FaHu,Hu,HuaZh,stripaper}.\par
For instance, Falconer and Hu, in \cite{FaHu}, considered Dirichlet problems defined on the Sierpi\'nski fractal.
 More precisely, under certain hypotheses on the nonlinear term, the existence of at least one non-trivial solution
was proved (see Theorems 3.5 and 3.18 of \cite{FaHu} and Remark \ref{Remark3} below).\par

 Further, in \cite{Hu}, Hu analyzed the following problem
  \begin{equation}\tag{$S_{a,\lambda}^{f}$}\label{NHu}
\left\{
\begin{array}{l}
\Delta u(x)+a(x)u(x)=\lambda f(x,u(x)),\quad x \in K\setminus K_0, \\
u|_{K_0}=0,\\
\end{array}
\right.
\end{equation}

\noindent  where $K$ is the Sierpi\'nski gasket (of intrinsic boundary $K_0$) in $\erre^2$ and $f:K\times\erre \to \erre$ is a continuous symmetric function satisfying some monotonicity properties. More precisely, in Theorem 2.2 of the above cited work, the existence of $p$--pairs of non-trivial solutions of \eqref{NHu} was achieved in relation with the value of the $p$--th eigenvalue, say $\lambda_p$, of the problem
$$
\left\{
\begin{array}{l}
\Delta u(x)+\lambda a(x)u(x)=0,\quad x \in K\setminus K_0, \\
u|_{K_0}=0.\\
\end{array}
\right.
$$
 \indent Very recently, Breckner, Repov\v{s} and Varga \cite{BRV} studied the existence of multiple solutions for the problem (\ref{N}) through variational methods. Their approach ensures the existence of at least three weak solutions under some hypotheses on the behaviour of the nonlinearity $f$.\par
 \indent Successively, Breckner, R\u{a}dulescu, and Varga \cite{BRaduV} proved the existence of infinitely many solutions of problem $(S_{a,1}^{f,g})$ under the key assumption that the nonlinearity $f$ is non-positive in a sequence of positive intervals (see Remark \ref{infinitely}). See also the papers \cite{GMB2, GMB3, GMB4, BrVa, BrVa2, MolRad} and references therein for related topics.\par
\indent In this paper, requiring an asymptotic behaviour of the nonlinearity $f$ at zero, we are able to determine a precise open interval of positive parameters $\lambda$, for
which problem \eqref{N} admits at least one non-trivial weak solution in the Sobolev space $H^1_0(V)$.\par
The proofs of our main results are based on a critical point theorem due to Ricceri \cite{ricceri} in the form given in \cite{GMB}.

 \begin{theorem}\label{BMB}
Let $X$ be a reflexive real Banach space, and let $\Phi,\Psi:X\to\R$
be two sequentially weakly lower semicontinuous functionals. Assume also that $\Phi$ is strongly continuous and coercive. For every $r>\inf_X
\Phi$, put
$$
\varphi(r):=\inf_{u\in\Phi^{-1}(]-\infty,r[)}\frac{\displaystyle\left(\sup_{v\in\Phi^{-1}(]-\infty,r[)}\Psi(v)\right)-\Psi(u)}{r-\Phi(u)}.
$$
Then, for every $r>\inf_X\Phi$ and every
$\lambda\in\left]0,{1}/{\varphi(r)}\right[$, the restriction
of the functional $I_\lambda:=\Phi-\lambda\Psi$ to
$\Phi^{-1}(]-\infty,r[)$ admits a global minimum, which is a
critical point $($local minimum$)$ of $I_\lambda$ in $X$.
\end{theorem}

{Clearly, the abstract framework introduced in the above mentioned paper is adaptable to our context by using the geometric and analytic properties of the Sierpi\'nski fractal as, for instance, the careful analysis of the Sobolev-type inequality (see, for instance, \cite[Lemma 2.4]{FaHu} and Section 2)
\begin{equation}\label{si}
\displaystyle\sup_{x,y\in V_*}\frac{|u(x)-u(y)|}{|x-y|^{\sigma}}\leq (2N+3)\sqrt{W(u)},
\end{equation}
\noindent where
$$
\sigma:=\displaystyle\frac{\log((N+2)/N)}{2\log 2},
$$
and $V_*$ and $W$ will be defined in the sequel.\par
}
\indent A special case of our results reads as follows.\par
\begin{theorem}\label{intro}
Let $f:\erre\to\erre$ be a non-negative continuous function. Assume that
$$
\lim_{\xi\rightarrow0^+} \frac{\di\int_0^\xi f(t)dt}{\xi^2}=+\infty.\eqno(\rm{h}_0')
$$
\indent Then the positive number $\lambda^{*}$, given by
$$
\lambda^{*}:=-\frac{1}{2(2N+3)^2\left(\di\int_Vg(x)d\mu\right)}
\sup_{\gamma>0}\frac{\gamma^2}{\displaystyle \int_0^\gamma f(t)dt},
$$
is such that, for every $\lambda\in ]0,\lambda^{*}[$,
the elliptic Dirichlet problem \eqref{N}
admits at least one non-trivial weak solution $u_{\lambda}\in (H^{1}_0(V),\|\cdot\|)$. Furthermore, $\displaystyle\lim_{\lambda\rightarrow 0^+}\|u_\lambda\|=0$.
\end{theorem}

 This paper is organized as follows. In Section 2 we recall the geometrical construction of the Sierpi\'nski gasket and our variational framework. Successively, Section 3 is devoted to the proof of the main theorem. Finally, in the last section, we give an application of the obtained results.\par
 We cite the very recent monograph by Krist\'aly, R\u adulescu and Varga \cite{KRV} as a general reference for the basic notions used here.

\section{Abstract Framework}

Let $N\geq2$ be a natural number and
let $p_1,\dots, p_N\in\R^{N-1}$ be so that $|p_i-p_j|=1$ for $i\neq j$. Define, for every
$i\in\{1,\dots,N\}$, the map $S_i\colon\R^{N-1}\to\R^{N-1}$ by
$$S_i(x)=\frac12\,x+\frac12\,p_i\,.$$
Let
${\mathcal S}:=\{S_1,\dots, S_N\}$ and denote by $F:{
\mathbb{P}}(\R^{N-1})\to{\mathbb{P}}(\R^{N-1})$ the map assigning to a subset
$A$ of $\R^{N-1}$ the set
$$F(A)=\bigcup_{i=1}^NS_i(A).$$
It is known that there
is a unique non-empty compact subset $V$ of $\R^{N-1}$, called the
\textit{attractor} of the family ${\mathcal S}$, such that $F(V)=V$; see, Theorem 9.1 in Falconer \cite{Fa}.\par
 The set $V$ is called
the \textit{Sierpi\'nski gasket} in $\R^{N-1}$ of \textit{intrinsic boundary} $V_0:=\{p_1,\dots,p_N\}$. Let $\mu$ be the normalized restriction of the $d$-dimensional Hausdorff measure ${\mathcal
H}^d$ on
$\R^{N-1}$ to the subsets of $V$, so $\mu(V)=1$.\par
 Further, the following property
of $\mu$ will be useful in the sequel:
\begin{equation}\label{support}
\mu(B)>0,\hbox{ for every non-empty open subset $B$ of }V.
\end{equation}
\indent In other words, the support of $\mu$ coincides with $V$ (see, for instance, Breckner, R\u adulescu and Varga \cite{BRaduV} for more details).\par
\indent Denote by
$C(V)$ the space of real-valued continuous functions on $V$ and by
$$C_0(V):=\{u\in C(V)\mid u|_{V_0}=0\}.$$
The spaces $C(V)$ and $C_0(V)$ are endowed with the usual supremum
norm $\|\cdot\|_\infty$. For a function $u\colon V\to\R$ and for $m\in\N$
let
\begin{equation}\label{defWm}
W_m(u)=\left(\frac{N+2}{N}\right)^m\sum_{\underset{|x-y|=2^{-m}}{x,y\in
V_m}}(u(x)-u(y))^2,
\end{equation}
where $V_m:=F(V_{m-1})$, for $m\geq1$. Put $V_*:=\bigcup_{m\geq 0} V_m$ and note that $V=\overline{V_*}$.\par
We have $W_m(u)\leq W_{m+1}(u)$ for very natural $m$, so we can put
\begin{equation}\label{defW}
W(u)=\lim_{m\to\infty} W_m(u).
\end{equation}
Define now
$$H_0^1(V):=\{u\in C_0(V)\mid W(u)<\infty\}.$$

It turns out that $H_0^1(V)$ is a dense linear subset of
$L^2(V,\mu)$ equipped with the $\|\cdot\|_2$ norm. We now endow
$H_0^1(V)$ with the norm
$$\|u\|:=\sqrt{W(u)}.$$
In fact, there is an inner product defining this norm: for $u,v\in
H_0^1(V)$ and $m\in\N$ let
$${\mathcal W}_m(u,v)=\left(\frac{N+2}{N}\right)^m\sum_{\underset{|x-y|=2^{-m}}{x,y\in V_m}}(u(x)-u(y))(v(x)-v(y)).$$
Put
$${\mathcal W}(u,v)=\lim_{m\to\infty} {\mathcal W}_m(u,v).$$

\noindent Then ${\mathcal W}(u,v)\in\R$ and the space $H_0^1(V)$, equipped with the
inner product ${\mathcal W}$, which induces the norm $\|\cdot\|$,
becomes a real Hilbert space.\par \indent Moreover,
\begin{equation}\label{embeddingconstant}
\|u\|_\infty\leq(2N+3)\|u\|, \hbox{ for every }u\in H_0^1(V),
\end{equation}
and the embedding
\begin{equation}\label{embedding}
(H_0^1(V),\|\cdot\|)\hookrightarrow (C_0(V),\|\cdot\|_\infty)
\end{equation}
is compact.\par
For
more details concerning the definitions and notions which lead in a natural way
to the Sobolev space $H^1_0(V)$ we refer to Fukushima and Shima \cite{FuSc}. See also
Sections 1.3 and 1.4 of \cite{Str2} (this reference applies for $N = 3$, but the cases
$N\geq 4$ are straightforward generalizations of this one).

\begin{remark}\label{equivalent}
\rm{
As pointed out by Falconer and Hu \cite{FaHu}, we just observe that if $a\in L^1(V)$ and $a\leq 0$ in $V$ then, by (\ref{embeddingconstant}),
the norm $\|\cdot\|_*$, defined by
$$
\|u\|_*:=\left({\mathcal W}(u,u)-\int_Va(x)u(x)^2d\mu\right)^{1/2},\quad \forall u\in H^1_0(V),
$$
is equivalent to $\|\cdot\|$.
}
\end{remark}

\indent We now state a useful property of the space $H_0^1(V)$ which shows,
together with the facts that $(H_0^1(V),\|\cdot\|)$ is a Hilbert
space and  $H_0^1(V)$ is dense in $L^2(V,\mu)$, that ${\mathcal
W}$ is a Dirichlet form on $L^2(V,\mu)$.

\begin{lemma}\label{MaMianalogue}
Let $h\colon \R\to\R$ be a Lipschitz mapping with Lipschitz constant
$L\geq0$ and such that $h(0)=0$. Then, for every $u\in H_0^1(V)$, we
have $h\circ u\in H_0^1(V)$ and $\|h\circ u\|\leq L\|u\|$.
\end{lemma}
\begin{pf}
It is clear that $h\circ u\in C_0(V)$. For every $m\in\N$ we have,
by (\ref{defWm}) and the Lipschitz property of $h$, that
$$W_m(h\circ u)\leq L^2 W_m(u).$$
Hence $W(h\circ u)\leq L^2 W(u)$, according to (\ref{defW}). Thus
$h\circ u\in H_0^1(V)$ and $\|h\circ u\|\leq L\|u\|$.
\end{pf}

\indent Following Falconer and Hu \cite{FaHu},
we can define in a standard way a linear self-adjoint operator $\Delta\colon Z\to L^2(V,\mu)$,
where $Z$ is a linear subset of $H_0^1(V)$ which is dense in $L^2(V,\mu)$ (and dense also in $(H_0^1(V),\|\cdot\|)$),
such that
$$-{\mathcal W}(u,v)=\int_V\Delta u\cdot vd\mu,\hbox{ for every }(u,v)\in Z\times H_0^1(V).$$
The operator $\Delta$ is called the {\textit{$($weak $)$ Laplacian} on}
$V$.\par \indent More precisely, let $H^{-1}(V)$ be the closure of
$L^2(V,\mu)$ with respect to the pre-norm
$$
\|u\|_{-1}=\sup_{\underset{\|h\|=1}{h\in H^1_0(V)}} |<u,h>|,
$$
\noindent where
$$
<v,h>=\int_Vv(x)h(x)d\mu,
$$
\noindent $v\in L^2(V,\mu)$ and $h\in H^1_0(V)$. Then $H^{-1}(V)$ is a Hilbert space. Then the relation
$$
-{\mathcal W}(u,v)=<\Delta u,v>,\,\,\,\,\,\,\forall v\in H^1_0(V),
$$
uniquely defines a function $\Delta u\in H^{-1}(V)$ for every $u\in
H^1_0(V)$.\par \indent Finally, fix $\lambda>0$. Let $a\colon
V\to\R$, $f\colon\R\to\R$ and $g\colon V\to\R$ be as in the
Introduction. We say that a function $u\in H_0^1(V)$ is a
{\textit{weak solution} of} (\ref{N}) if
$${\mathcal W}(u,v)-\int_Va(x)u(x)v(x)d\mu+\lambda\int_Vg(x)f(u(x))v(x)d\mu=0,$$
for every $v\in H_0^1(V)$.\par

While we mainly work with the weak Laplacian, there is also a directly
defined version. We say that $\Delta_su$ is the \textit{standard Laplacian} of $u$ if $\Delta_su:V\to\R$ is continuous and
$$
\lim_{m\to\infty}\sup_{x\in V\setminus V_0}|(N+2)^m(H_mu)(x)-\Delta_su(x)|=0,
$$
\noindent where
$$
(H_mu)(x):=\sum_{\underset{|x-y|=2^{-m}}{y\in V_m}}(u(y)-u(x)),
$$
\noindent for $x\in V_m$. We say that $u\in C_0(V)$ is a {\it strong solution} of (\ref{N}) if $\Delta_su$ exists and is continuous for all $x\in V\setminus V_0,$ and
$$
\Delta_s u(x)+a(x)u(x)=\lambda g(x)f(u(x)),\,\,\,\,\,\forall\; x\in V\setminus V_0.
$$
\indent
The existence of the standard Laplacian of a function $u\in H^1_0(V)$ implies
the existence of the weak Laplacian $\Delta u$ (see, for completeness, Falconer and Hu \cite{FaHu}).\par
\begin{remark}\label{strong}
\rm{
If $a\in C(V)$, $f\colon\R\to\R$ is continuous and $g\in C(V)$,
then,  using  the regularity result Lemma 2.16 of Falconer and Hu \cite{FaHu}, it
follows that every weak solution of the problem (\ref{N}) is also a
strong solution.}
\end{remark}
\section{Main results}

 Define $F\colon\R\to\R$ by $F(\xi)=\displaystyle\int_0^\xi f(t)dt$ and fix $\lambda>0$.
The functional $I_\lambda\colon H_0^1(V)\to\R$  given by
\begin{equation}\label{energyfunctional}
I_\lambda(u):=\frac12\|u\|^2-\frac12\int_Va(x)u(x)^2d\mu+\lambda\int_Vg(x)F(u(x))d\mu,\
\end{equation}
for every $u\in H_0^1(V)$,
will turn out to be the energy functional attached to problem
(\ref{N}).\par
We have the following result contained in \cite[Proposition 2.19]{FaHu} that we recall here in a convenient form.
\begin{lemma}\label{criticalpoint}
The energy functional $I_\lambda\colon H_0^1(V)\to\R$ defined by relation
{\rm (\ref{energyfunctional})} is a $C^1(H_0^1(V),\R)$ functional.
Moreover, for each point $u\in H^1_0(V)$,
$$I_\lambda'(u)(v)={\mathcal W}(u,v)-\int_Va(x)u(x)v(x)d\mu+\lambda\int_Vg(x)f(u(x))v(x)d\mu,\,\,\,\forall\; v\in H^1_0(V).$$
 In particular, $u\in H_0^1(V)$ is a weak solution of problem \eqref{N}
if and only if $u$ is a critical point of $I_\lambda$.
\end{lemma}

\indent The aim of the paper is to prove the following result concerning the
existence of at least one non-trivial solutions of the problem (\ref{N}).

\begin{theorem}\label{Main1}
Let $f:\erre\to\erre$ be a continuous function with $f(0)=0$. Assume that

$$
\displaystyle -\infty<\liminf_{\xi\rightarrow 0^+}\frac{F(\xi)}{\xi^{2}}\,\,\,\,\,\,\,\ \textrm{and}\,\,\,\,\,\,\,\,\,\limsup_{\xi\rightarrow 0^{+}
}\frac{F(\xi)}{\xi^{2}}=+\infty.\eqno{(\rm{h}_{0})}
$$

\noindent Then the positive number $\lambda^{*}$, given by
$$
\lambda^{*}:=-\frac{1}{2(2N+3)^2\left(\di\int_Vg(x)d\mu\right)}
\sup_{\gamma>0}\frac{\gamma^2}{\displaystyle\max_{|\xi|\leq {\gamma}}F(\xi)},
$$
is such that, for every $\lambda\in ]0,\lambda^{*}[$,
the problem \eqref{N} admits at least one non-trivial weak solution $u_{\lambda}\in H_0^1(V)$. Moreover,
$$
\lim_{\lambda\rightarrow 0^+}\|u_{\lambda}\|=0,
$$
and, for every $\bar \gamma>0$, the function $\lambda\rightarrow I_{\lambda}(u_\lambda)$ is negative and strictly decreasing in
$$\left]0,
-\frac{1}{2(2N+3)^2\left(\di\int_Vg(x)d\mu\right)}
\frac{\bar\gamma^2}{\displaystyle\max_{|\xi|\leq {\bar\gamma}}F(\xi)}\right[.
$$
\end{theorem}

\begin{pf}
Let us define the functionals
$\Phi,\Psi:X\rightarrow\RR$ by
$$\Phi(u):=\frac12\|u\|^2-\frac12\int_Va(x)u(x)^2d\mu \ \ \ \mbox{and}\ \ \ \Psi(u):=-\int_Vg(x)F(u(x))d\mu,$$
 where $X$ denotes the reflexive Banach space $H^1_0(V)$. Now, in order to achieve our goal, fix $\lambda$ as in the conclusion.\par
 \indent With the above notations we have that $I_\lambda=\Phi-\lambda\Psi$.
We seek for
weak solutions of problem (\ref{N}) by applying Theorem \ref{BMB}. First of all we observe that, by Lemma \ref{criticalpoint}, the functional $I_\lambda\in C^{1}(X,\R)$.\par
 Moreover, $\Phi$ is obviously coercive and, by using Lemma 5.6 in Breckner, R\u adulescu and Varga \cite{BRaduV}, the functionals $\Phi$ and $\Psi$ are weakly sequentially lower
semicontinuous on $X$.\par
 Since $0<\lambda<\lambda^{*}$,
there exists $\bar{\gamma}>0$ such that
\begin{equation}\label{n3}
\lambda<\lambda^{*}{(\bar{\gamma})}:=-\frac{\bar{\gamma}^2}{2(2N+3)^2\left(\di\int_Vg(x)d\mu\right)\displaystyle\max_{|\xi|\leq \bar{\gamma}}F(\xi)}.
\end{equation}
\indent
Set $\displaystyle r:=\frac{\bar{\gamma}^2}{2(2N+3)^2}$.
Due to the compact embedding into $C_0(V)$, by (\ref{embeddingconstant}), we have
$$
\left\{v\in X: \Phi(v)<r \right\}\subseteq \left\{v\in X: \|v\|_\infty\leq \bar{\gamma} \right\}.
$$
Therefore
\begin{eqnarray*}
\varphi(r) &=& \inf_{\Phi(u)<r}\frac{\displaystyle\sup_{\Phi(v)<
r}\displaystyle\int_V(-g(x))F(v(x))d\mu+\int_Vg(x)F(u(x))d\mu}{r-\Phi(u)}\nonumber\\
 &\leq&\frac{\displaystyle\sup_{\Phi(v)<
r}\displaystyle\int_V(-g(x))F(v(x))d\mu}{r}\nonumber\\
&\leq&-\left(\int_Vg(x)d\mu\right)\frac{\displaystyle\max_{|\xi|\leq \bar{\gamma}}F(\xi)}{r}\nonumber\\
&=&-2(2N+3)^2\left(\di\int_Vg(x)d\mu\right)\frac{\displaystyle\max_{|\xi|\leq \bar{\gamma}}F(\xi)}{\bar{\gamma}^{2}}=\frac{1}{\lambda^{*}{(\bar{\gamma})}}.\nonumber
\end{eqnarray*}

\indent Thanks to Theorem \ref{BMB}, there exists a function $u_\lambda\in\Phi^{-1}(]-\infty,r[)$ such that
$$
I'_{\lambda}(u_\lambda)=\Phi'(u_\lambda)-\lambda\Psi'(u_\lambda)=0,
$$
and, in particular, $u_\lambda$ is a global minimum of the restriction of $I_{\lambda}$ to $\Phi^{-1}(]-\infty,r[)$. \indent Now, we claim that the function $u_\lambda$ cannot be trivial, i.e. $u_\lambda\neq 0$. Indeed, fix a non-negative function $u\in X$ such that there is an element $x_0\in
V$ with $u(x_0)>1$. It follows that
 $$
 D:=\{x\in V\mid u(x)>1\}
 $$
 is a non-empty open subset of $V$ (due to the continuity of $u$).\par
  Define $h\colon\R\to\R$ as follows
$$
h(t):=|\min\{t,1\}|,\qquad\mbox{for all}\ t\in\R.
$$
Then $h(0)=0$ and $h$ is a Lipschitz function whose Lipschitz constant $L$ is equal to $1$. Hence, by using Lemma
\ref{MaMianalogue}, it follows that $v:=h\circ u\in X$.\par
 Moreover,
$v(x)=1$ for every $x\in D$, and $0\leq v(x)\leq1$ for every $x\in
V$.\par
\noindent On the other hand, condition
$$
-\infty <\liminf_{\xi\rightarrow 0^+}\frac{F(\xi)}{\xi^2}
$$
 implies the existence of
real numbers $\rho>0$ and $\varrho$ such that
\begin{equation}\label{mic4}
F(\xi)\geq \varrho\xi^2,\hbox{ for every } \xi\in[0,\rho[.
\end{equation}

Further, condition
$$
\limsup_{\xi\rightarrow 0^+
}\frac{F(\xi)}{\xi^{2}}=+\infty
$$
 yields the existence of a sequence $\{\xi_n\}$
in $]0,\rho[$ such that $\displaystyle\lim_{n\to\infty}\xi_n=0$ and
\begin{equation}\label{mic5}
\lim_{n\to\infty}\frac{F(\xi_n)}{\xi_n^2}=+\infty.
\end{equation}

Now, we have that
\begin{eqnarray*}
I_\lambda(\xi_nv)=\frac{\xi_n^2}{2}\|v\|^2-\frac{\xi_n^2}{2}\int_Va(x)v(x)^2d\mu&+&\lambda F(\xi_n)\int_Dg(x)d\mu\\
&+&\lambda \int_{V\setminus D}g(x)F(\xi_nv(x))d\mu,
\end{eqnarray*}
for every $n\in\N$.\par

\indent Using (\ref{mic4}) and the fact that $g\leq0$ in $V$, we get
\begin{eqnarray*}
I_\lambda(\xi_nv)\leq \frac{\xi_n^2}{2}\|v\|^2-\frac{\xi_n^2}{2}\int_Va(x)v(x)^2d\mu&+&\lambda F(\xi_n)\int_Dg(x)d\mu\\
&+&\lambda\varrho\xi_n^2\int_{V\setminus D}g(x)v(x)^2d\mu,
\end{eqnarray*}
for
every $n\in\N$. Thus
$$\frac{I_\lambda(\xi_nv)}{\xi_n^2}\leq \frac12\|v\|^2-\frac12\int_Va(x)v(x)^2d\mu+\lambda\frac{F(\xi_n)}{\xi_n^2}\int_Dg(x)d\mu+\lambda\varrho\int_{V\setminus D}g(x)v(x)^2d\mu.$$

Condition $({\rm{h}}_2)$ and (\ref{support}) imply that $$\int_Dg(x)d\mu<0,$$
so we get from (\ref{mic5}) and the above inequality that
$$\lim_{n\to\infty}\frac{I_\lambda(\xi_nv)}{\xi_n^2}=-\infty.$$

Then, there is an index $n_0$ such that $I_\lambda(\xi_nv)<0$ for every $n\geq
n_0$. Now, since
$$
\lim_{n\rightarrow \infty}\Phi(\xi_nv)=0,
$$
one has that $\xi_nv\in \Phi^{-1}(]-\infty,r[)$ definitively. In conclusion, $0_{X}$ cannot be a global minimum for the restriction of the functional $I_\lambda$ to $\Phi^{-1}(]-\infty,r[)$. Hence, for every $\lambda\in ]0,\lambda^{*}[$ the problem \eqref{N} admits a non-trivial solution $u_{\lambda}\in X$.\par
 At this point, we prove that $\|u_{\lambda}\|\rightarrow 0$ as $\lambda\rightarrow 0^+$ and that
the function $\lambda\rightarrow I_{\lambda}(u_\lambda)$ is negative and decreasing in $]0,\lambda^{*}[$.\par
For our goal, let us consider $\bar{\lambda}\in ]0,\lambda^{*}[$. Moreover, let $\bar{\gamma}>0$ and let $\lambda\in ]0,\lambda^{*}(\bar{\gamma})[$. The functional $I_{\lambda}$ admits a non-trivial critical point $u_{\lambda}\in \Phi^{-1}(]-\infty,r[)$, where $$\displaystyle r:=\frac{\bar{\gamma}^2}{2(2N+3)^2}.$$

 Since $\Phi$ is coercive and $u_\lambda\in \Phi^{-1}(]-\infty,r[)$ for every
$\lambda\in ]0,\lambda^{*}(\bar{\gamma})[$, there exists a positive number $L$ such that
$$
\|u_\lambda\|\leq L,
$$
for every
$\lambda\in ]0,\lambda^{*}(\bar{\gamma})[$.\par
\indent Therefore, since $\Psi'$ is a compact operator, there exists a positive constant $M$ such that
\begin{equation}\label{bounded}
\left|\Psi(u_{\lambda})\right|\leq \|\Psi'(u_{\lambda})\|_{X^*}\|u_{\lambda}\|< ML^2,
\end{equation}
\noindent for every $\lambda\in ]0,\lambda^{*}(\bar{\gamma})[$.\par

\indent Now $I'_\lambda(u_{\lambda})=0$, for every $\lambda\in ]0,\lambda^{*}(\bar{\gamma})[$ and in particular $$I'_\lambda(u_{\lambda})(u_{\lambda})=0,$$ that is,
\begin{equation}\label{9}
\Phi(u_\lambda)=\lambda\int_V g(x)f(u_{\lambda}(x))u_{\lambda}(x)d\mu,
\end{equation}
\noindent for every $\lambda\in ]0,\lambda^{*}(\bar{\gamma})[$.\par
 Hence, by (\ref{bounded}) and (\ref{9}) it follows that
\begin{equation}\label{c}
\lim_{\lambda\rightarrow 0^+}\Phi(u_\lambda)=0.
\end{equation}
Moreover, one has
\begin{equation}\label{a}
\frac{\|u_\lambda\|^2}{2}\leq \frac{\|u_\lambda\|^2}{2}- \frac{\displaystyle\int_Va(x)u_\lambda(x)^2d\mu}{2}=\Phi(u_\lambda),
\end{equation}
\noindent for every $\lambda\in ]0,\lambda^{*}(\bar{\gamma})[$. Then, conditions (\ref{c}) and (\ref{a}) yield
$$
\lim_{\lambda\rightarrow 0^+}{\|u_{\lambda}\|}=0.$$

\indent Further, the map $\lambda\mapsto I_{\lambda}(u_\lambda)$ is negative in $]0,\lambda^{*}(\bar{\gamma})[$ since the restriction
of the functional $I_{\lambda}$ to
$\Phi^{-1}(]-\infty,r[)$ admits a global minimum, which is a
critical point $($local minimum$)$ of $I_{\lambda}$ in $X$.\par

\indent Finally, observe that
$$
I_{\lambda}(u)=\lambda\left(\frac{\Phi(u)}{\lambda}-\Psi(u)\right),
$$
for every $u\in X$ and fix $0<\lambda_1<\lambda_2<\lambda^{*}(\bar{\gamma})$.\par
 Moreover, put
$$
m_{\lambda_1}:=\left(\frac{\Phi(u_{\lambda_1})}{\lambda_1}-\Psi(u_{\lambda_1})\right)=\inf_{u\in\Phi^{-1}(]-\infty,r[)}\left(\frac{\Phi(u)}{\lambda_1}-\Psi(u)\right),
$$
and
$$
m_{\lambda_2}:=\left(\frac{\Phi(u_{\lambda_2})}{\lambda_2}-\Psi(u_{\lambda_2})\right)=\inf_{u\in\Phi^{-1}(]-\infty,r[)}\left(\frac{\Phi(u)}{\lambda_2}-\Psi(u)\right).
$$
Clearly, as claimed before, $m_{\lambda_i}<0$ (for $i=1,2$), and $m_{\lambda_2}\leq m_{\lambda_1}$ thanks to $\lambda_1<\lambda_2$.

Then the map $\lambda\mapsto I_{\lambda}(u_{\lambda})$ is strictly decreasing in $]0,\lambda^{*}(\bar{\gamma})[$ owing to
$$
I_{\lambda_2}(u_{\lambda_2})=\lambda_2m_{\lambda_2}\leq \lambda_2m_{\lambda_1}<\lambda_1m_{\lambda_1}=I_{\lambda_1}(u_{\lambda_1}).
$$
The proof is complete.
\end{pf}

\begin{remark} \label{rem2}
\rm{
\noindent We observe that condition ($\rm{h}_0$) is technical and ensures that the solution, obtained by using Theorem \ref{Main1}, is non-trivial. Anyway, the statements of Theorem \ref{Main1} are still true for every continuous function $f$ that does not vanish at zero. In this last case our approach ensures the existence of one non-trivial solution, for $\lambda\in ]0,\lambda^{*}[$, without condition ($\rm{h}_0$). If
$$
\displaystyle\max_{|\xi|\leq {\bar\gamma}}F(\xi)=0,
$$
for some $\bar\gamma>0$, Theorem \ref{Main1} ensures the existence of one non-trivial solution, for every $\lambda\in ]0,+\infty[$.
}
\end{remark}

\begin{remark} \label{infinitely}
\rm{If in addition to condition ($\rm{h}_0$) in Theorem \ref{Main1}, the function $f$ also satisfies
\begin{itemize}
\item[($\rm{h}_1'$)] There exist two sequences $\{a_n\}$ and $\{b_n\}$ in $]0,\infty[$ with $b_{n+1}<a_n<b_{n}$,
$\displaystyle\lim_{n\to \infty}b_n=0$ and such that $f(s)\leq 0
\hbox{ for every }s\in [a_n,b_n];$
\item[($\rm{h}_2'$)] Either $\sup\{s<0\mid f(s)>0\}=0$, or there is a
$\delta>0$ with $f|_{[-\delta,0]}= 0$,
\end{itemize}
\noindent then, as proved by Breckner, R\u{a}dulescu and Varga in \cite{BRaduV}, the problem $(S_{a,1}^{f,g})$
admits a sequence $\{u_n\}$ of pairwise distinct weak solutions
such that
$\displaystyle\lim_{n\to\infty}\|u_n\|=0.$ In particular,
$\displaystyle\lim_{n\to\infty}\|u_n\|_\infty=0$.
}
\end{remark}
\begin{remark} \label{rem3}
\rm{
\noindent A sufficient condition that ensures hypothesis ($\rm{h}_0$) in Theorem \ref{Main1} is expressed by
$$
\lim_{\xi\rightarrow0^+} \frac{F(\xi)}{\xi^2}=+\infty.\eqno(\rm{h}_0')
$$
Further, if $f$ is non-negative, one has
$$
\sup_{\gamma>0}\frac{\gamma^2}{\displaystyle\max_{|\xi|\leq {\gamma}}F(\xi)}=\sup_{\gamma>0}\frac{\gamma^2}{\displaystyle F(\gamma)},
$$
since, in this case, $\displaystyle\max_{|\xi|\leq {\gamma}}F(\xi)=F(\gamma)$, for every positive $\gamma$.
Hence, Theorem \ref{intro} in the Introduction immediately follows from Theorem \ref{Main1}.
}
\end{remark}
The following example is a direct consequence of Theorem \ref{Main1}, bearing in mind Remarks \ref{rem2} and \ref{strong}.
\begin{example}
\rm{
\noindent For each
parameter $\lambda$ belonging to
$$
\Lambda:=\left]0,\frac{2e^{-2}}{(2N+3)^2}\right[,
$$
\noindent the following Dirichlet problem
\begin{equation}\tag{$S_{\lambda}$}\label{N3}
\left\{
\begin{array}{l}
\Delta u(x)+\lambda e^{u(x)}=u(x),\quad x \in V\setminus V_0, \\
u|_{V_0}=0,\\
\end{array}
\right.
\end{equation}
\noindent admits at least one non-trivial strong solution. Moreover,
$$
\lim_{\lambda\rightarrow 0^+}\|u_{\lambda}\|=0
$$
and the function $\lambda\rightarrow I_{\lambda}(u_\lambda)$ is negative and decreasing in $\Lambda$.}
\end{example}
\begin{remark}\label{Remark3}
\rm{
In \cite{FaHu} Falconer and Hu studied the non-autonomous Dirichlet problem
\begin{equation}\tag{$S_{a,\lambda}^{f}$}\label{N34}
\left\{
\begin{array}{l}
\Delta u(x)+a(x)u(x)=\lambda f(x,u(x)),\quad x \in V\setminus V_0, \\
u|_{V_0}=0,\\
\end{array}
\right.
\end{equation}
\noindent where $a:V\to\erre$ is assumed to be integrable and $f:V\times \erre\to\erre$ is a continuous function. The celebrated Ambrosetti-Rabinowitz condition
\begin{itemize}
\item[(AR)]\textit{there are constants $\nu>2$ and $r\geq 0$ such that}
$$
tf(x,t)\leq \nu F(x,t)<0,
$$
\textit{for every} $|t|\geq r$, {\it uniformly for every $x\in V$,}
\end{itemize}
 is an essential request in almost all the existence theorems contained in the above cited paper. However
if, for instance, $f$ is constant for large $|t|$, assumption (AR) is
violated, even though (\ref{N34}) would be expected to have a non-trivial
solution. The saddle point theorem copes with this case; see \cite[Theorem 4.2]{FaHu}.
We observe that Theorem \ref{Main1} (see also Remark \ref{rem2}) obtained in this paper does not require a global growth of the non-linearity $f$ in order to obtain the existence of one non-trivial solution as the above example shows.}
\end{remark}
{\bf Acknowledgements.}  This paper was written when G.M.B. was visiting professor at the University of Ljubljana in 2014. He expresses his gratitude to the host institution for warm hospitality.
   The manuscript was realized within the auspices of the INdAM - GNAMPA Project 2015 titled {\it Modelli ed equazioni non-locali di tipo frazionario} and the SRA grants P1-0292-0101 and J1-5435-0101.

\end{document}